\documentclass[preprint]{elsarticle}

\usepackage{float}

\usepackage{lineno,hyperref}
\modulolinenumbers[5]

\usepackage{amsmath}
\usepackage{amsfonts}
\allowdisplaybreaks[4]

\usepackage{tikz}
\usetikzlibrary{shapes.geometric, arrows}
\tikzstyle{process} = [rectangle, rounded corners, minimum width=4cm, minimum height=0.7cm, text centered, font=\small, draw=black, fill=white]
\tikzstyle{arrow} = [thick,->,>=stealth]

\usepackage{url}

\usepackage{multirow}
\usepackage{threeparttable}

\usepackage{mdwlist}

\usepackage{pgfplots}
\pgfplotsset{width=6cm,compat=1.13}

\journal{Renewable Energy}

\newcommand{\lj}[1]{\textcolor{black}{{}#1}}
\begin{document}

\begin{frontmatter}

\title{Impacts and Benefits of UPFC to Wind Power Integration in Unit Commitment}
% \tnotetext[mytitlenote]{Fully documented templates are available in the elsarticle package on \href{http://www.ctan.org/tex-archive/macros/latex/contrib/elsarticle}{CTAN}.}

%% Group authors per affiliation:
% \author{Elsevier\fnref{myfootnote}}
% \address{Radarweg 29, Amsterdam}
% \fntext[myfootnote]{Since 1880.}

%% or include affiliations in footnotes:
\author[mymainaddress]{Jia Li}
% \ead[url]{www.elsevier.com}

\author[mymainaddress]{Feng Liu\corref{mycorrespondingauthor}}
\cortext[mycorrespondingauthor]{Corresponding author}
\ead{lfeng@tsinghua.edu.cn}

\author[mysecondaryaddress]{Zuyi Li}

\author[mymainaddress]{Shengwei Mei}

\author[mytertiaryaddress]{Guangyu He}

\address[mymainaddress]{State Key Laboratory of Power Systems, Department of Electrical Engineering, Tsinghua University, Beijing 100084, China}
\address[mysecondaryaddress]{Robert W. Galvin Center for Electricity Innovation at Illinois Institute of Technology, Chicago, IL 60616 USA}
\address[mytertiaryaddress]{School of Electronic Information and Electrical Engineering, Shanghai Jiao Tong University, Shanghai 200240, China}

\begin{abstract}
Unified Power Flow Controller (UPFC) is recognized as the most powerful flexible AC transmission systems (FACTS) device for power system operation. This paper addresses how  UPFC explores the transmission flexibility and facilitates the integration of uncertain and volatile wind power generation. To this end, a comprehensive unit commitment (UC) model with UPFC and uncertain wind power generation is proposed. Then, some metrics are introduced to evaluate the impacts of UPFC on the reliability, security and economy of power system operation. Further, different dispatch strategies of UPFC are compared to provide helpful guidances on making full use of UPFC to hedge against uncertainties. In addition, facing the challenging mixed-integer non-linear non-convex problems, approximate models are proposed to provide a starting point to solve the problems efficiently. All these models are easy to adapt to other types of FACTS devices. Illustrative numerical results are provided.
 
\end{abstract}


\end{frontmatter}

% \linenumbers

\section*{Nomenclature}

\subsection*{Indexes}
\begin{basedescript}{\desclabelstyle{\pushlabel}\desclabelwidth{6em}}
 	\item [$i,j,k$] Index of system buses.
 	\item [$s$] Index of wind power generation scenarios.
 	\item [$t$] Index of time periods.
\end{basedescript}

\subsection*{Sets}
\begin{basedescript}{\desclabelstyle{\pushlabel}\desclabelwidth{8em}}
	\item [$N/N^L/N^G/N^W$] Set of all buses/load buses/thermal unit buses/wind farm buses.
	\item [$S$]	Set of wind power generation scenarios.
	\item [$T$]	Set of time periods.
\end{basedescript}

\subsection*{Parameters}
\begin{basedescript}{\desclabelstyle{\pushlabel}\desclabelwidth{9em}}
	\item [$\alpha^{LS}/\alpha^{WC}$]	Price of load shedding/wind power curtailment.
	\item [$B_{ij}/G_{ij}/G_{i}^s$]	Susceptance/conductance/ shunt conductance.
	\item [$p_s$]	Probability of scenario $s$.
	\item [$P_{ij}^{dc,\max}$]	Limitation of active power transferred through the converters of UPFC.
	\item [$P_{i,t}^L/Q_{i,t}^L$]	Active/reactive load.
	\item [$P_{i}^{G,\max}/Q_{i}^{G,\max}$]	Maximum active/reactive power output of unit $i$.
	\item [$P_{i}^{G,\min}/Q_{i}^{G,\min}$]	Minimum active/reactive power output of unit $i$.
	\item [$P_{ij}^{TC}$]	Capacity of line $ij$.
	\item [$P_{i,t}^{WF}/P_{s,i,t}^W$]	Active wind power generation forecast/scenarios.
	\item [$Q_{i,t}^{WF}/Q_{s,i,t}^W$]	Reactive load forecast/scenarios of wind farm.
	\item [$R_t$] Spinning reserve requirement.
	\item [$RD_i/RU_i$]	Ramp-down/ramp-up limit of unit $i$.
	\item [$SD_i/SU_i$]	Shutdown/startup ramp limit of unit $i$.
	% \item [$TN$]	Number of periods of the time span.
	\item [$T_{ij}^{se,\max}/T_{ij}^{sh,\max}$]	Thermal limitation of series/shunt converter.
	\item [$V^{\max}/V^{\min}$]	Maximum/minimum voltage magnitude.
	\item [$X_{ij}$] Reactance of line $ij$.
	\item [$\Delta_{ij}^{P,U}/\Delta_{ij}^{Q,se}/\Delta_{ij}^{Q,sh}$] Re-dispatch limit of active/reactive power injection of UPFC.
\end{basedescript}

\subsection*{Variables}
\begin{basedescript}{\desclabelstyle{\pushlabel}\desclabelwidth{6em}}
	\item [$C_{i,t}^D/C_{i,t}^U$]	Shutdown/startup cost of a thermal unit.
	\item [$C_{s,i,t}^F$]	Fuel cost of a thermal unit in the second stage.
	\item [$P_{i,t}^G/P_{s,i,t}^G$]	Active power generation of a thermal unit in the first/second stage.
	\item [$P_{i,t}^{G,avl}/P_{s,i,t}^{G,avl}$]	Maximum available active power generation of a thermal unit in the first/second stage.
	\item [$P_{ij}/P_{s,ij}$] Branch power flow in the first/second stage.
	\item [$P_{ij,t}^U/P_{s,ij,t}^U$] Active power injection of UPFC in the first/ second stage.
	\item [$P_{s,i,t}^{LS}$] Load shedding in the second stage.
	\item [$P_{s,i,t}^{WC}$] Wind power curtailment in the second stage.
	\item [$Q_{i,t}^G/Q_{s,i,t}^G$]	Reactive power generation of a thermal unit in the first/second stage.
	\item [$Q_{ij,t}^{U,sh}/Q_{s,ij,t}^{U,sh}$]	Shunt reactive power injection of UPFC in the first/second stage.
	\item [$Q_{ij,t}^{U,se}/Q_{s,ij,t}^{U,se}$]	Series reactive power injection of UPFC in the first/second stage.
	\item [$Q_{ij,t}^U/Q_{s,ij,t}^U$] Equivalent non-control reactive power injection of UPFC in the first/second stage.
	\item [$Q_{s,i,t}^{WC}$] Reactive load shedding of wind farm.
	\item [$u_{i,t}$] Thermal unit status.
	\item [$V_{i,t}/V_{s,i,t}$]	Voltage magnitude in the first/second stage.
	\item [$\theta_{ij,t}/\theta_{s,ij,t}$]	Voltage angle in the first/second stage.
\end{basedescript}

\section{Introduction}\label{sec_introduction}

The uncertainty of wind power generation has posed new challenges to power systems. The inherent volatility of wind power generation may impact the security and economy of power system operation, causing voltage violation and congestion. In order to deal with the increasing penetration of wind power generation, there is an urgent need to take full advantage of power system flexibility. 

Generally, the flexibility of power system operation can be divided into three categories: generation side, transmission network and demand side. In generation side and demand side, a lot of efforts have been devoted to addressing the uncertainty of wind power generation. Different kinds of methods have been applied in unit commitment (UC) and economic dispatch to enhance generation-side flexibility, such as stochastic optimization \cite{Wang2008,Zhang2011}, chance-constrained optimization \cite{Wang2012}, robust optimization \cite{Jiang2012a,Ye2015c}, minimax regret \cite{Jiang2013a}. In order to improve the computational efficiency and tractability, a scenario tree approach \cite{Nasri2015} and a scenario selection algorithm inspired by importance sampling \cite{PapavasiliouA.andOren2013,Papavasiliou2012a,Papavasiliou2010,Papavasiliou2012} have been applied to characterize the uncertainty of wind power generation with a few scenarios in stochastic UC problems. A practical adaptive robust UC solution methodology has been proposed in \cite{Bertsimas2013}, highlighting the scalability of the proposed formulation. In demand side, demand response, including price-based and incentive-based methods has been introduced to improve demand-side flexibility. 

Although flexible AC transmission systems (FACTS) and high-voltage direct current (HVDC) provide transmission flexibility, few studies have been made to analyze their impacts on wind power integration in power system operation. In terms of FACTS, a control scheme is proposed in \cite{Yang2012} to determine the optimal steady-state settings of thyristor controlled series capacitor (TCSC) in order to improve the usage of the existing transmission network. A two-stage method based on regression analysis is applied using a collection of offline simulations. A scenario-based optimal power flow (OPF) model is proposed in \cite{Nasri2014} to minimize wind power spillage with TCSC. The decision making method is formulated as a two-stage stochastic programming model while the control of TCSC are considered in the second stage. It has been proved that a series compensation may reduce the wind power spillage, unserved load, and total active power losses of the network \cite{Nasri2014}. In \cite{Thakurta2015}, nodal and angular sensitivities are used for coordinating phase shifting transformers (PSTs) to deal with contingencies and to increase wind power penetration. It has been shown that the coordination in operation of these devices helps to bring the system into a secure state from an overload situation. The flexibility of HVDC lines and HVDC grid is exploited in security-constrained OPF frameworks so as to minimize the operating cost under wind power uncertainty \cite{Vrakopoulou2013a,Vrakopoulou2013}. In those papers, the problems are formulated as chance constrained optimization programs and scenario-based methodologies are applied, which offer a strong solution with a-priori constraint violation guarantees. It has been shown that HVDC lines can be used to handle the fluctuating in-feed from renewable energy sources \cite{Vrakopoulou2013}. In \cite{Rabiee2014}, a stochastic multi-period OPF model is presented which consists of an offshore wind farm connected to the grid by a line-commutated converter HVDC link. The obtained results demonstrate that the availability of transmission network capacity at the interface of AC/DC network is a key factor affecting the utilization of wind power generation.

\lj{Though efforts have been devoted to introducing FACTS devices into economic dispatch \cite{ISI:000387207300010,Mukherjee2016,Shchetinin2016,ISI:000387881500004}, the references with regard to unit commitment with FACTS devices are still limited. Reference \cite{Sreejith2015} focuses on solving security constrained UC problem using artificial bee colony (ABC) algorithm incorporating FACTS devices. The results show that the installation of FACTS devices can improve power flow and reduce transmission line losses. A UC model considering FACTS devices for corrective operation is proposed in \cite{Sahraei-Ardakani2015a}. The original mixed-integer non-linear problem is reformulated as an mixed-integer linear problem. In this context, though optimality is not guaranteed, the simulation studies show that the method finds the optimal solution in most cases. However, these references all focus on using FACTS devices to improve power system reliability considering contingency, other than using FACTS devices to promote the integration of wind power generation. In addition, no evaluation process is developed for comparing different ways of using FACTS devices. In \cite{Sreejith2015}, different types of FACTS devices are modeled in detail. In this paper, though only one type of FACTS devices, i.e., UPFC is considered, the power injection model is applied. It is more general and can easily be used to model other types of FACTS devices.
In \cite{Sahraei-Ardakani2015a}, a linear programming approach is proposed to reduce computational complexity. However, only DC power flow is considered. In this context, the flexibility of FACTS device may not be fully exploited. In this paper, both active and reactive power flow are taken into account in order to take full advantage of FACTS devices.}

Compared with generation-side and demand-side flexibility, using FACTS devices may be faster and cheaper, making it an appropriate tool to cope with the uncertainty of wind power generation in UC. FACTS devices can be initiated quickly and frequently, since the power electronics allows very short reaction time down to far below one second \cite{Zhang2012a}. Moreover, FACTS devices are capable of controlling the interrelated parameters that govern the operation of transmission systems including series impedance, shunt impedance, current, voltage, phase angle \cite{Hingorani2000a}. Among the converter-based FACTS devices, Unified Power Flow Controller (UPFC) \cite{Gyugyi1995,Sen1998} is a versatile FACTS device, which is capable of controlling circuit impedance, voltage angle and power flow simultaneously for optimal operation performance of power system \cite{Zhang2012a}. 

\lj{This paper is focused on analyzing the impacts of UPFC on wind power integration in UC, so as to provide new insights into the utilization of UPFC. In the context, two issues need to be addressed. First, how to assess the impacts of UPFC? Second, how to use UPFC in an appropriate way? To lay the foundation for evaluation, a comprehensive UC model with UPFC and uncertain wind power generation is proposed. Then, some metrics are introduced to evaluate the impacts of UPFC. Further, different dispatch strategies of UPFC are compared to facilitate wind power integration. Additionally, facing the challenging mixed-integer non-linear non-convex problems, approximate models are proposed to provide a starting point to solve the problems efficiently. Thus, the contributions of this paper are threefold:
} 

\lj{(1) A comprehensive evaluation process based on a two-stage stochastic UC model with UPFC and uncertain wind power generation is proposed. The control variables of UPFC are incorporated in both stages, while AC power flows are taken into account. In addition, power injection model of UPFC is used, which can be extended easily to incorporate different types of FACTS devices, such as static var compensator (SVC), static synchronous compensator (STATCOM), thyristor-controlled phase-shifter (TCPS).
}

\lj{(2) Different ways of applying UPFC in UC are comprehensively compared to fully exploit the flexibility of UPFC. Numerical results show that using UPFC only in the first stage helps to make a more economic UC schedule, but may bring adverse effects on wind power integration. On the other hand, the expectation of wind power curtailment and load shedding are largely reduced with the help of UPFC in the second stage. Meanwhile, the voltage profile is improved and the expected total cost is reduced even within a more rigid voltage limit. Moreover, the lowest expected total cost is achieved when UPFC is dispatched in both stages. The results provide new insights for the use of UPFC in UC so as to address uncertain wind power generation.}

\lj{(3) A DC model and a mixed model are proposed to find an approximate solution with lower computational burdens. Since the proposed AC model is a mixed-integer non-linear non-convex problem which is challenging to solve, a good starting point is identified to solve it. The DC model is computationally efficient but it cannot make full use of UPFC. Meanwhile, the mixed model obtains relatively accurate solutions with less computational burdens than the AC model, when the UC schedule derived from the DC model is feasible.
}

The rest of the paper is organized as below. Section \ref{sec_upfc_uc} incorporates UPFC into the two-stage stochastic UC model, and proposes the approximate DC and mixed models. Section \ref{sec_decision} proposes the different strategies of using UPFC. Section \ref{sec_evaluation} presents the evaluation process and introduces various metrics. Section \ref{sec_case} provides the comprehensive analysis of the impacts of UPFC. Finally, in Section \ref{sec_conclusion}, some relevant conclusions are drawn.

\section{Two-stage UC Model with UPFC}\label{sec_upfc_uc}
\subsection{Two-Stage Stochastic UC Model}
The UC problem is formulated as a two-stage stochastic programming model. The UC schedule is determined in the first stage, while wind power curtailment and load shedding are only allowed in the second stage.

The objective is to minimize the total cost, which consists of two parts: the UC cost in the first stage, including startup cost and shutdown cost, and the expected cost in the second stage, including fuel cost, wind power curtailment cost and load shedding cost (\ref{eq_obj}).

\begin{equation}
\begin{aligned}
	\min & \sum_{t\in T}\sum_{i\in N^G}\left(C_{i,t}^U + C_{i,t}^D\right) \\
	& + \sum_{s\in S}p_s\left(C_{s,i,t}^F + \sum_{t\in T}\sum_{i\in N^W}\alpha^{WC}P_{s,i,t}^{WC}\right. \\
	& \left. + \sum_{t\in T}\sum_{i\in N^L}\alpha^{LS}P_{i,t}^{LS}\right) \label{eq_obj}
\end{aligned}
\end{equation}

\subsubsection{First-stage Problem}
The first-stage problem represents the here-and-now decision making process before knowing the actual values of stochastic variables. Therefore, the decisions are made based on forecast data. The constraints are formulated with reference to \cite{Carrion2006}, including active and reactive power balance constraints (\ref{eq_first_p_balance})-(\ref{eq_first_q_balance}), voltage magnitude limits (\ref{eq_first_v}), voltage angle limits (\ref{eq_first_theta}), spinning reserve requirements (\ref{eq_first_reserve}), transfer capacity limits (\ref{eq_first_p_line})-(\ref{eq_first_line_cap}), thermal unit generation limits (\ref{eq_first_pg_min})-(\ref{eq_first_qg}), ramping constraints (\ref{eq_first_ramp_up})-(\ref{eq_first_ramp_down}), minimum up and down time constraints \cite{Carrion2006}.

\begin{align}
	P_{i,t}^G & - P_{i,t}^L + P_{i,t}^{WF} \nonumber \\ 
	& = V_{i,t}\sum_{j\in N}V_{j,t}\left(G_{ij}\cos\theta_{ij,t} + B_{ij}\sin\theta_{ij,t}\right), \nonumber \\
	& \forall i\in N, \forall t\in T \label{eq_first_p_balance} \\
	Q_{i,t}^G & - Q_{i,t}^L + Q_{i,t}^{WF} \nonumber \\ 
	& = V_{i,t}\sum_{j\in N}V_{j,t}\left(G_{ij}\sin\theta_{ij,t} - B_{ij}\cos\theta_{ij,t}\right), \nonumber \\
	& \forall i\in N, \forall t\in T \label{eq_first_q_balance} \\
	& V^{\min} \le V_{i,t} \le V^{\max}, \forall i\in N, \forall t\in T \label{eq_first_v} \\
	& -\pi \le \theta_{ij,t} \le \pi, \forall i,j\in N, \forall t\in T \label{eq_first_theta} \\
	& \sum_{i\in N^G}P_{i,t}^{G,avl} + \sum_{t\in N^W}P_{i,t}^{WF} \ge \sum_{i\in N}P_{i,t}^L + R_t, \forall t\in T \label{eq_first_reserve} \\
	P_{ij} & = V_{i,t}^2(-G_{ij} + G_i^s) + V_{i,t}V_{j,t}G_{ij}\cos\theta_{ij,t} \nonumber \\
	& + V_{i,t}V_{j,t}B_{ij}\sin\theta_{ij,t}, \forall i,j\in N, \forall t\in T \label{eq_first_p_line} \\
	& |P_{ij}| \le P_{ij}^{TC}, \forall i,j \in N \label{eq_first_line_cap} \\
	& P_{i}^{G,\min}u_{i,t} \le P_{i,t}^G \le P_{i,t}^{G,avl}, \forall i\in N^G, \forall t\in T \label{eq_first_pg_min} \\
	& 0 \le P_{i,t}^{G,avl} \le P_{i}^{G,\max}u_{i,t}, \forall i\in N^G, \forall t\in T \\
	& Q_{i}^{G,\min}u_{i,t} \le Q_{i,t}^G \le Q_{i}^{G,\max}u_{i,t}, \nonumber \\
	& \forall i\in N^G, \forall t\in T \label{eq_first_qg} \\
	P_{i,t}^{G,avl} & \le P_{i,t-1}^G + RU_i u_{i,t-1} + SU_i(u_{i,t} - u_{i,t-1}) \nonumber \label{eq_first_ramp_up} \\
	& + P_i^{G,\max}(1-u_{i,t}), \forall i\in N^G, \forall t\in T \\
	P_{i,t}^{G,avl} & \le P_{i}^{G,\max}u_{i,t+1} + SD_i(u_{i,t} - u_{i,t+1}), \nonumber \\
	& \forall i\in N^G, \forall t\in T \\
	P_{i,t-1}^G & \le P_{i,t}^G + RD_i u_{i,t} + SD_i(u_{i,t-1} - u_{i,t}) \nonumber \\
	& + P_i^{G,\max}(1-u_{i,t-1}), \forall i\in N^G, \forall t\in T \label{eq_first_ramp_down}
\end{align}

\subsubsection{Second-stage Problem} \label{sec_second_stage}
The second-stage problem makes the wait-and-see decisions, based on wind power generation scenarios. Latin hypercube sampling (LHS) technique \cite{glasserman2003monte} is employed to generate a set of scenarios in order to form a discrete approximation of wind power generation. However, a large number of scenarios may make the associated stochastic optimization problem intractable. Thus, a scenario reduction technique based on probability metric \cite{Morales2009,Growe-Kuska2003} is applied to reduce the number of scenarios, while preserving most of the stochastic information.

In each scenario, the power system operation constraints are satisfied, including active and reactive power balance constraints (\ref{eq_second_p_balance})-(\ref{eq_second_q_balance}), wind power curtailment limits (\ref{eq_wc}), load shedding limits (\ref{eq_ls}), voltage magnitude limits (\ref{eq_second_v}), voltage angle limits (\ref{eq_second_theta}), spinning reserve requirements (\ref{eq_second_reserve}), transfer capacity limits (\ref{eq_second_p_line})-(\ref{eq_second_line_cap}), thermal unit generation limits (\ref{eq_second_pgmin})-(\ref{eq_second_qg}), and ramping constraints (\ref{eq_second_ramp_up})-(\ref{eq_second_ramp_down}).

\begin{align}
	P_{s,i,t}^G & - P_{s,i,t}^L + P_{s,i,t}^{LS} + P_{s,i,t}^{W} - P_{s,i,t}^{WC} & \nonumber \\ 
	& = V_{s,i,t}\sum_{j\in N}V_{s,j,t}\left(G_{ij}\cos\theta_{s,ij,t} + B_{ij}\sin\theta_{s,ij,t}\right), \nonumber \\
	& \forall s\in S, \forall i\in N, \forall t\in T \label{eq_second_p_balance} \\
	Q_{s,i,t}^G & - Q_{s,i,t}^L - Q_{s,i,t}^{W} + Q_{s,i,t}^{WC} & \nonumber \\ 
	& = V_{s,i,t}\sum_{j\in N}V_{s,j,t}\left(G_{ij}\sin\theta_{s,ij,t} - B_{ij}\cos\theta_{s,ij,t}\right), \nonumber \\
	& \forall s\in S, \forall i\in N, \forall t\in T \label{eq_second_q_balance} \\
	& 0 \le P_{s,i,t}^{WC} \le P_{s,i,t}^{W}, \forall s\in S, \forall i\in N^W, \forall t\in T \label{eq_wc} \\
	& 0 \le P_{s,i,t}^{LS} \le P_{i,t}^L, \forall s\in S, \forall i\in NL, \forall t\in T \label{eq_ls} \\
	& V^{\min}\le V_{s,i,t} \le V^{\max}, \nonumber \\
	& \forall s\in S, \forall i\in N, \forall t\in T \label{eq_second_v} \\
	& -\pi \le \theta_{s,ij,t} \le \pi,  \forall s\in S, \forall i,j\in N, \forall t\in T \label{eq_second_theta} \\
	& \sum_{i\in N^G}P_{s,i,t}^{G,avl} + \sum_{t\in N^W}(P_{s,i,t}^{W} - P_{s,i,t}^{WC}) \nonumber \\
	& \ge \sum_{i\in N}(P_{i,t}^L - P_{s,i,t}^{LS}) + R_t, \nonumber \\
	&  \forall s\in S, \forall t\in T \label{eq_second_reserve} \\
	P_{s,ij} & = V_{s,i,t}^2(-G_{ij} + G_i^s) + V_{s,i,t}V_{s,j,t}G_{ij}\cos\theta_{s,ij,t} \nonumber \\
	& + V_{s,i,t}V_{s,j,t}B_{ij}\sin\theta_{s,ij,t}, \forall s\in S, \nonumber \\
	& \forall i,j\in N, \forall t\in T \label{eq_second_p_line} \\
	& |P_{s,ij}| \le P_{ij}^{TC}, \forall i,j \in N \label{eq_second_line_cap} \\
	& P_{i}^{G,\min}u_{i,t} \le P_{s,i,t}^G \le P_{s,i,t}^{G,avl}, \nonumber \\
	& \forall s\in S, \forall i\in N^G, \forall t\in T \label{eq_second_pgmin} \\
	& 0 \le P_{s,i,t}^{G,avl} \le P_{i}^{G,\max}u_{i,t}, \nonumber \\
	& \forall s\in S, \forall i\in N^G, \forall t\in T \\
	& Q_{i}^{G,\min}u_{i,t} \le Q_{s,i,t}^G \le Q_{i}^{G,\max}u_{i,t}, \nonumber \\
	& \forall s\in S, \forall i\in N^G, \forall t\in T \label{eq_second_qg} \\
	P_{s,i,t}^{G,avl} & \le P_{s,i,t-1}^G + RU_i u_{i,t-1} + SU_i(u_{i,t} - u_{i,t-1}) \nonumber \\
	& + P_i^{G,\max}(1-u_{i,t}), \forall s\in S, \forall i\in N^G, \forall t\in T \label{eq_second_ramp_up} \\
	P_{s,i,t}^{G,avl} & \le P_{i}^{G,\max}u_{i,t+1} + SD_i(u_{i,t} - u_{i,t+1}), \nonumber \\
	& \forall s\in S, \forall i\in N^G, \forall t\in T \\
	P_{s,i,t-1}^G  & \le P_{s,i,t}^G + RD_i u_{i,t} + SD_i(u_{i,t-1} - u_{i,t}) \nonumber \\
	& + P_i^{G,\max}(1-u_{i,t-1}), \forall s\in S, \forall i\in N^G, \forall t\in T \label{eq_second_ramp_down}
\end{align}

\subsection{Incorporation of UPFC in UC}
\begin{figure}[htb]
\centering
\includegraphics[width=2in]{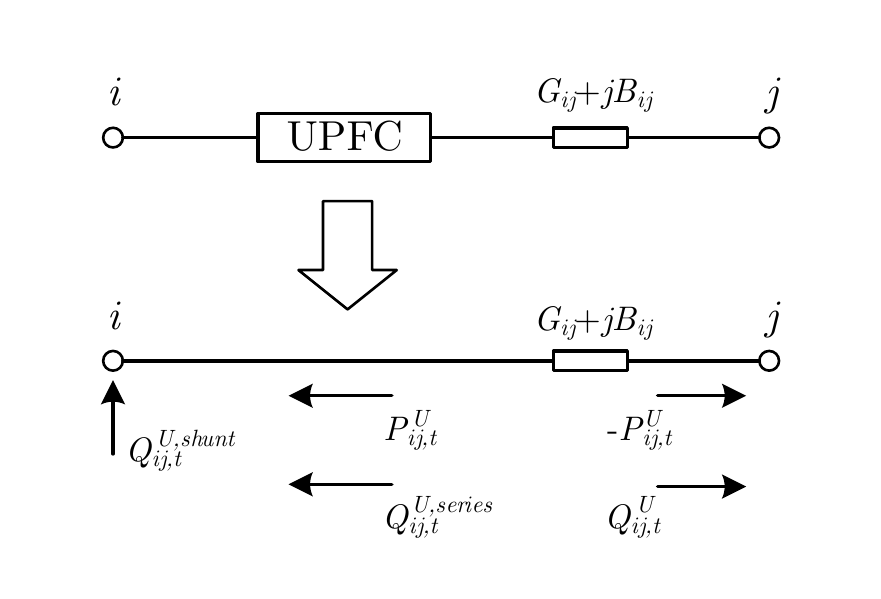}
\caption{Power injection model of UPFC.}
\label{fig_upfc}
\end{figure}
The steady-state models of FACTS devices can be formulated as a series and/or shunt-inserted voltage (-current) source(s), which is called voltage source model (VSM) \cite{song1999flexible}. The VSM is intuitive but it destroys the symmetric characteristics of the admittance matrix \cite{Han1982}. Derived from VSM, the power injection model (PIM) \cite{Noroozian1997,Song2002}, as shown in Fig. \ref{fig_upfc}, keeps the symmetry of the admittance matrix. The PIM of UPFC is applied in both first and second stages. Both active and reactive power flows are taken into account, so that not only the impacts on active power but also voltage can be analyzed. For brevity, \lj{only the integration of UPFC into the first stage of UC is shown}. 
UPFC is included in the second stage analogously. Then, re-dispatch constraints are introduced between the two stages to formulate different dispatch strategies.

Let $L$ be the set of branches with UPFC installed. For any branch $ij\in L$, the PIM of UPFC in the first stage is formulated as below \cite{Song2002}.

\begin{align}
	P_{i,t}^G & - P_{i,t}^L + P_{i,t}^{WF} + P_{ij,t}^U \nonumber \\
	& = V_{i,t}\sum_{k\in N}V_{k,t}\left(G_{ik}\cos\theta_{ik,t} + B_{ik}\sin\theta_{ik,t}\right), \forall t\in T \label{eq_upfc_pi} \\
	Q_{i,t}^G & - Q_{i,t}^L + Q_{i,t}^{WF} + Q_{ij,t}^{U,se} + Q_{ij,t}^{U,sh} \nonumber \\
	& = V_{i,t}\sum_{k\in N}V_{k,t}\left(G_{ik}\sin\theta_{ik,t} - B_{ik}\cos\theta_{ik,t}\right), \forall t\in T \label{eq_upfc_qi} \\
	P_{j,t}^G & - P_{j,t}^L + P_{j,t}^{WF} + P_{ij,t}^U \nonumber \\
	& = V_{j,t}\sum_{k\in N}V_{k,t}\left(G_{jk}\cos\theta_{jk,t} + B_{jk}\sin\theta_{jk,t}\right), \forall t\in T \label{eq_upfc_pj} \\
	Q_{j,t}^G & - Q_{j,t}^L + Q_{j,t}^{WF} + Q_{j,t}^{U} \nonumber \\
	& = V_{j,t}\sum_{k\in N}V_{k,t}\left(G_{jk}\sin\theta_{jk,t} - B_{jk}\cos\theta_{jk,t}\right), \forall t\in T \label{eq_upfc_qj} \\
	P_{ij} & = V_{i,t}^2\left(-G_{ij} + G_i^s\right) + V_{i,t}V_{j,t}G_{ij}\cos\theta_{ij,t} \nonumber \\
	& + V_{i,t}V_{j,t}B_{ij}\sin\theta_{ij,t} - P_{ij,t}^U, \forall t\in T \label{eq_upfc_p_line} \\
	& \sqrt{\left(Q_{ij,t}^{U,sh}\right)^2 + \left(P_{ij,t}^U\right)^2} \le T_{ij}^{sh,\max}, \forall t\in T \label{eq_upfc_t_shunt} \\
	& \sqrt{\left(Q_{ij,t}^{U,se}\right)^2 + \left(P_{ij,t}^U\right)^2} \le T_{ij}^{se,\max}, \forall t\in T \\
	& |P_{ij,t}^U| \le P_{ij}^{dc,\max}, \forall t\in T \label{eq_upfc_p_dc}
\end{align}

The impacts of UPFC on active power flow are incorporated into the active power balance constraint (\ref{eq_first_p_balance}) as two inverse active power injections at bus $i$ (\ref{eq_upfc_pi}) and bus $j$ (\ref{eq_upfc_pj}), respectively. Similarly, The impacts of UPFC on reactive power flow are included in the reactive power balance constraint (\ref{eq_first_q_balance}) as reactive power injections at bus $i$ (\ref{eq_upfc_qi}) and bus $j$ (\ref{eq_upfc_qj}), respectively. It should be noted that the two reactive power injections have different physical meanings. According to \cite{Song2002}, $Q_{ij,t}^{U,sh}$ is injected by the shunt synchronous voltage source (SVS) directly into bus $i$ for regulating voltage, while $Q_{ij,t}^{U,se}$ is generated by the series SVS, flowing via line $ij$ for reactive line flow control. Moreover, the formulation with two reactive power injections is more versatile, which can be extended to represent shunt or series controllers. Equality (\ref{eq_upfc_p_line}) represents the line flow with UPFC power injections. Inequality (\ref{eq_upfc_t_shunt})-(\ref{eq_upfc_p_dc}) denote the thermal limitations and the limit of active power transferred through converters. The details of UPFC injection model can be found in \cite{Noroozian1997}.

The model of UPFC in the second stage is formulated similarly to the above constraints (\ref{eq_upfc_pi})-(\ref{eq_upfc_p_dc}) with second-stage variables. In addition, the re-dispatch constraints are added as below, so that the power injections of UPFC are dispatched within an acceptable level between the two stages to accommodate wind power uncertainty.

\begin{align}
	& \left|P_{s,ij,t}^U - P_{ij,t}^U\right| \le \Delta_{ij}^{P,U}, \forall s\in S, \forall t\in T, \forall ij\in L \label{eq_upfc_d_p} \\
	& \left|Q_{s,ij,t}^{U,se} - Q_{ij,t}^{U,se}\right| \le \Delta_{ij}^{Q,se}, \forall s\in S, \forall t\in T, \forall ij\in L \\
	& \left|Q_{s,ij,t}^{U,sh} - Q_{ij,t}^{U,sh}\right| \le \Delta_{ij}^{Q,sh}, \forall s\in S, \forall t\in T, \forall ij\in L \label{eq_upfc_d_q_shunt}
\end{align}

\subsection{DC and Mixed Models}
The proposed AC problem above is mixed-integer, nonlinear, and non-convex. Due to its complexity and lack of efficient computational tools, \lj{some approximate models are proposed}, which may be served as alternatives to solve the problem. If only DC power flows are considered, the model can be simplified to reduce computational burdens. However, it is worth mentioning that the impacts of UPFC on voltage are ignored, yielding the flexibility of UPFC cannot be fully exploited. Compared with the AC model, the objective function is the same as (\ref{eq_obj}), while the constraints related to voltage and reactive power are removed in the DC model. Moreover, the power balance constraints and line flow equations in the first (\ref{eq_dc_first_p_balance})-(\ref{eq_dc_first_p_line_upfc}) and second stage (\ref{eq_dc_second_p_balance})-(\ref{eq_dc_second_p_line_upfc}) are changed, while the other constraints remain unchanged.

\begin{align}
	& P_{i,t}^G - P_{i,t}^L + P_{i,t}^{WF} = \sum_{j\in N}\frac{\theta_{ij,t}}{X_{ij}}, \forall t\in T, \forall ij\notin L \label{eq_dc_first_p_balance} \\
	& P_{ij} = \frac{\theta_{ij,t}}{X_{ij}}, \forall t\in T, \forall ij\notin L \\
	& P_{i,t}^G - P_{i,t}^L + P_{i,t}^{WF} + P_{ij,t}^U = \sum_{k\in N}\frac{\theta_{ik,t}}{X_{ik}}, \forall t\in T, \forall ij\in L \\
	& P_{j,t}^G - P_{j,t}^L + P_{j,t}^{WF} - P_{ij,t}^U = \sum_{k\in N}\frac{\theta_{jk,t}}{X_{jk}}, \forall t\in T, \forall ij\in L \\
	& P_{ij} = \frac{\theta_{ij,t}}{X_{ij}} - P_{ij,t}^U, \forall t\in T, \forall ij\in L \label{eq_dc_first_p_line_upfc} \\
	& P_{s,i,t}^G - P_{s,i,t}^L + P_{s,i,t}^{LS} + P_{s,i,t}^{W} - P_{s,i,t}^{WC} = \sum_{j\in N}\frac{\theta_{s,ij,t}}{X_{ij}}, \nonumber \label{eq_dc_second_p_balance} \\
	& \forall s\in S, \forall t\in T, \forall ij\notin L \\
	& P_{s,ij} = \frac{\theta_{s,ij,t}}{X_{ij}}, \forall s\in S, \forall t\in T, \forall ij\notin L \\
	& P_{s,i,t}^G - P_{i,t}^L + P_{s,i,t}^{LS} + P_{s,i,t}^{W} - P_{s,i,t}^{WC} + P_{s,ij,t}^U = \sum_{k\in N}\frac{\theta_{s,ik,t}}{X_{ik}}, \nonumber \\
	& \forall s\in S, \forall t\in T, \forall ij\in L \\
	& P_{s,j,t}^G - P_{j,t}^L + P_{s,i,t}^{LS} + P_{s,i,t}^{W} - P_{s,i,t}^{WC} - P_{s,ij,t}^U = \sum_{k\in N}\frac{\theta_{s,jk,t}}{X_{jk}}, \nonumber \\
	& \forall s\in S, \forall t\in T, \forall ij\in L \\
	& P_{s,ij} = \frac{\theta_{s,ij,t}}{X_{ij}} - P_{s,ij,t}^U, \forall s\in S, \forall t\in T, \forall ij\in L \label{eq_dc_second_p_line_upfc}
\end{align}

Since the integer variables largely increase computational complexity, the mixed model performs an AC economic dispatch based on the UC solution of the aforementioned DC model. First, the DC problem is solved to obtained the UC schedule. Second, the binary variables in the AC model, which indicate the unit status, are fixed with the UC schedule of the DC model. Third, the AC problem is solved with economic dispatch constraints. On one hand, the computational burden is reduced compared with the full AC model. On the other hand, the solution to the DC-only model may violate the AC economic dispatch constraints, since reactive power flow constraints are ignored in the DC model. 

\section{Different Dispatch Strategies of UPFC}\label{sec_decision}

In the first stage, the setpoints of thermal unit outputs and UPFC are determined according to wind power generation forecast minimizing the UC cost, while in the second-stage, the expected total cost is minimized with the re-dispatch of thermal units and UPFC to accommodate uncertain wind power generation in different scenarios. According to the stage where the control of UPFC is employed, different dispatch strategies are proposed as follows, aimed at seeking the best way to utilize UPFC for wind power integration in UC. All the models are implemented in GAMS \cite{GAMS} and solved by DICOPT \cite{Grossmann2001}.

\lj{(1) DM: Deterministic UC model without UPFC, based on wind power generation forecast.}

\lj{(2) NOM: Non-optimal model with UPFC, where only wind power generation forecast is used, and no optimization is made for the second stage.}

(3) NM: No UPFC model.

(4) FSM: UPFC controllable in the first stage model.

(5) SSM: UPFC controllable in the second stage model.

(6) FSSM: UPFC controllable in the first stage and second stage model.

\subsection{No UPFC Model (NM)}
In NM, all of UPFC associated variables are set to zero as below. The objective function and other constraints remain unchanged. As a result, it is a basic two-stage stochastic UC model with uncertain wind power generation, serving as a benchmark.

\begin{align}
	& P_{ij,t}^U = Q_{ij,t}^{U,se} = Q_{ij,t}^{U,sh} = Q_{ij,t}^U = 0, \forall t\in T, \forall ij\in L \\
	& P_{s,ij,t}^U = Q_{s,ij,t}^{U,se} = Q_{s,ij,t}^{U,sh} = Q_{s,ij,t}^U = 0, \nonumber \\
	& \forall s\in S, \forall t\in T, \forall ij\in L
\end{align}

\subsection{UPFC in the First Stage Model (FSM)}
When employed only in the first stage, the UPFC cannot be re-dispatched in the second stage. In other words, the setpoints of UPFC are determined in the first stage and remain unchanged in the second stage. Thus, the re-dispatch constraints of UPFC are set to zero as below.

\begin{equation}
	\Delta_{ij}^{P,U} = \Delta_{ij}^{Q,se} = \Delta_{ij}^{Q,sh} = 0, \forall s\in S, \forall t\in T, \forall ij\in L
\end{equation}

\subsection{UPFC in the Second Stage Model (SSM)}
When employed only in the second stage, the UPFC has zero setpoints in the first stage as below, and is re-dispatched in the second stage with respect to wind power generation scenarios.

\begin{equation}
	P_{ij,t}^{U} = Q_{ij,t}^{U,se} = Q_{ij,t}^{U,sh} = Q_{ij,t}^U = 0, \forall t\in T, \forall ij\in L
\end{equation}

\subsection{UPFC in the First and Second Stage Model (FSSM)}
When used in the both stages, the UPFC is set up in the first stage and then re-dispatched in the second stage based on wind power generation scenarios, while the re-dispatch constraints (\ref{eq_upfc_d_p})-(\ref{eq_upfc_d_q_shunt}) are satisfied.

It is worth noting that compared with NM, FSM has extra controllable variables in the first stage which may help to reduce operating cost, while SSM has extra controllable variables in the second stage which may contribute to reducing the expected wind power curtailment cost and load shedding cost. Further, FSSM has the most controllable variables among the proposed models, making it the most flexible one. 

\section{Evaluation and Metrics}\label{sec_evaluation}
\begin{figure}[htb]
\centering
\begin{tikzpicture}[node distance=1cm]
\node (sg) [process] {Scenario Generation};
\node (sr) [process, below of=sg] {Scenario Reduction};
\node (opt) [process, below of=sr] {Optimization};
\node (eva) [process, below of=opt] {Evaluation};
\node (met) [process, below of=eva] {Metrics Calculation};
\draw [arrow] (sg) -- (sr);
\draw [arrow] (sr) -- (opt);
\draw [arrow] (opt) -- (eva);
\draw [arrow] (eva) -- (met);
\end{tikzpicture}
\caption{Whole analysis process}
\label{fig_process}
\end{figure}
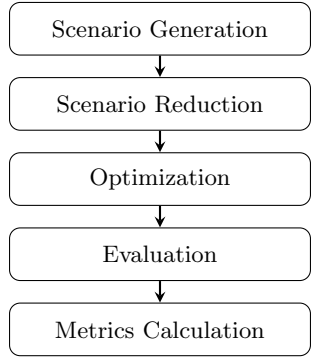

As stated in Section \ref{sec_second_stage}, a large quantity of possible wind power generation scenarios are generated using LHS, then reduced to a few scenarios by scenario reduction technique. The optimization problems are formulated based on the reduced scenarios. Therefore, evaluations are required to test the optimal solutions in each of the original generated scenario, as well as to analyze metrics reflecting the impacts of UPFC on wind power integration. In this paper, 1000 scenarios are generated then reduced to 10 for optimization. With the first-stage decisions fixed as the optimal solutions, including thermal unit status and UPFC setpoints, economic dispatch is performed according to each wind power generation scenario as a simulation of the second-stage decision making process aimed at minimizing operating cost, including fuel cost, wind power curtailment cost and load shedding cost. After all the scenarios are evaluated, various metrics are calculated based on all the evaluation results. The whole process is shown in Fig. \ref{fig_process}. 

The expected costs are calculated as below, including expected fuel cost (EFC), expected wind power curtailment cost (EWC), expected load shedding cost (ELC) and expected total cost (ETC).

\begin{align}
	\text{EFC} &= \sum_{s\in \hat{S}} p_s C_{s}^F \\
	\text{EWC} &= \sum_{s\in \hat{S}} p_s C_s^{WC} \\
	\text{ELC} &= \sum_{s\in \hat{S}} p_s C_s^{LS} \\
	\text{ETC} &= \text{EFC} + \text{EWC} + \text{ELC} + \text{UCC}
\end{align} 
where $\hat{S}$ denotes the set of all generated scenarios. $C_s^F$, $C_s^{WC}$, $C_s^{LS}$ are the fuel cost, wind power curtailment cost, and load shedding cost of each scenario, respectively. UC cost (UCC) includes the startup and shuntdown cost of thermal units.

Change rate (CR) of expected costs compares the difference of the EFC, EWC, ELC, ETC before and after UPFC is employed. In other words, it shows the rate of change of the expected costs in FSM, SSM and FSSM, compared with NM. For instance, the change rate of EFC ($\text{CR}^{EFC}$) in FSSM is calculated as below.

\begin{equation}
	\text{CR}^{EFC} = \frac{\text{EFC}^{FSSM} - \text{EFC}^{NM}}{\text{EFC}^{NM}}
\end{equation}
where the superscript of EFC denotes the type of model. The change rate of EWC ($\text{CR}^{EWC}$), ELC ($\text{CR}^{ELC}$) and ETC ($\text{CR}^{ETC}$) can be obtained similarly.

The loss of load probability (LOLP) is also introduced to evaluate the probability of load shedding (\ref{eq_lolp}).

\begin{equation}
	\text{LOLP} = \sum_{s\in \hat{S}}p_s\sum_{t\in T}\frac{h_t^L}{24} \label{eq_lolp}
\end{equation}
where $h_t^L$ equals to 1 if there is load shedding at hour $t$, otherwise $h_t^L$ is 0. $24$ indicates the dispatch horizon is 24 hours. Similarly, \lj{the wind power curtailment probability (WPCP) is proposed} to quantify the probability of wind power curtailment (\ref{eq_wpcp}).

\begin{equation}
	\text{WPCP} = \sum_{s\in \hat{S}}p_s\sum_{t\in T}\frac{h_t^W}{24} \label{eq_wpcp}
\end{equation}
where $h_t^W$ equals to $1$ if there is wind power curtailment at hour $t$, otherwise $h_t^W$ is $0$.

\section{Case Studies}\label{sec_case}
\begin{figure}[htb]
\centering
\includegraphics[width=2in]{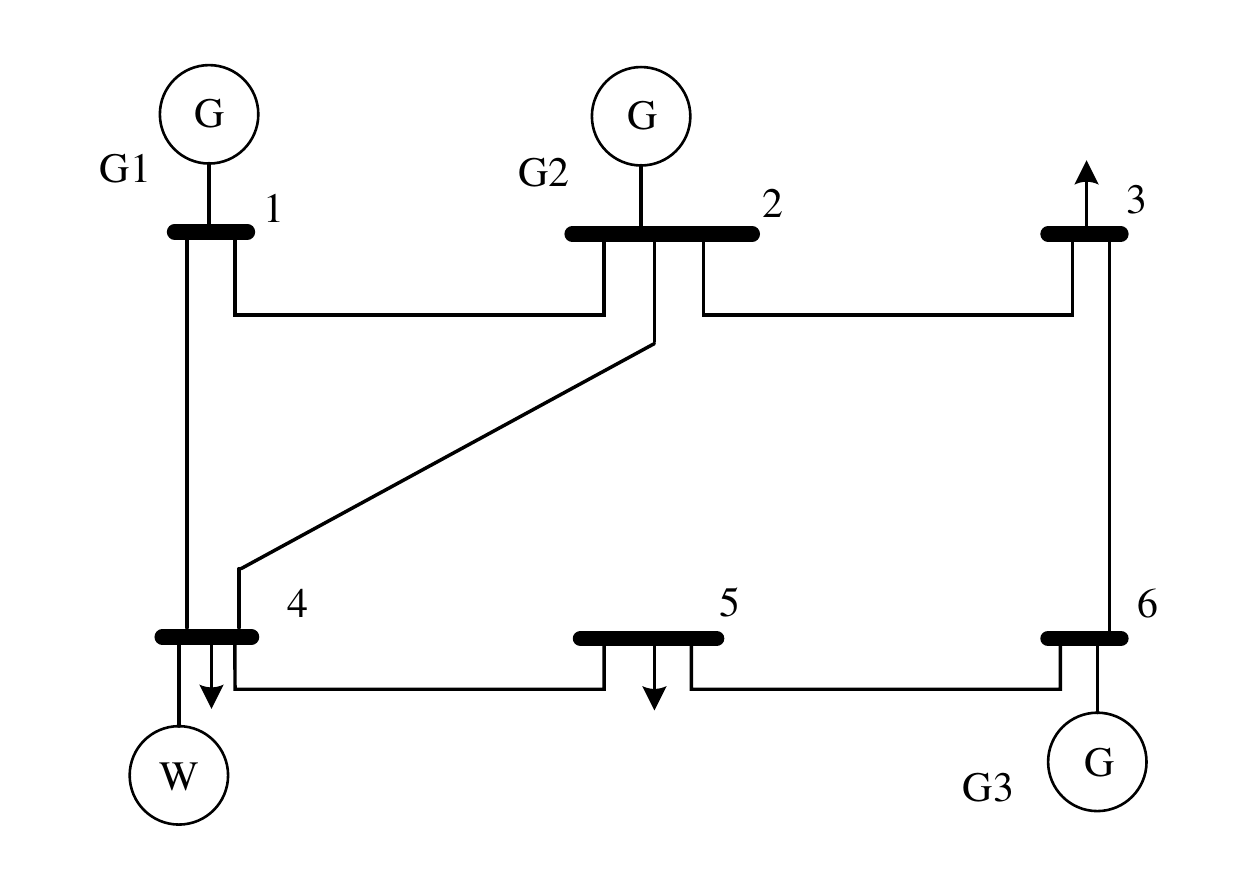}
\caption{Six-bus system.}
\label{fig_case6}
\end{figure}

A 6-bus system \cite{Wang2008}, as shown in Fig. \ref{fig_case6}, is used for testing the proposed models and analyzing the impacts of UPFC. The system contains three thermal units, one wind farm and one UPFC. The transmission line data are listed in Table \ref{tab_line}. The spinning reserve requirements are assumed to be 5\% of the load. The wind farm, which is assumed to be controlled with a constant power factor of 0.96, is located at bus 4, with the capacity of 150MW. The UPFC is installed in line 4-5 and paralleled at bus 4. The parameters of UPFC are given in Table \ref{tab_upfc}. The price of wind power curtailment is considered as the levelized cost of electricity \cite{U.S.EnergyInformationAdministration}, which is \$73.6/MWh, and the price of load shedding is assumed to be \$300/MWh. In Appendix, the characteristics of thermal units are given in Table \ref{tab_gen_data} and Table \ref{tab_operation_data}, while the hourly load and wind power generation forecast are listed in Table \ref{tab_load}.

\begin{table}[htbp]\footnotesize
  \centering
  \caption{Transmission line data.}
    \begin{tabular}{ccccccc}
    \hline
    Line & From & To & R & X & b & Capacity \\
    No. & Bus & Bus & (p.u.) & (p.u.) & (p.u.) & (MW) \\
    \hline
    1     & 1     & 2     & 0.005 & 0.170  & 0.02  & 150 \\
    2     & 1     & 4     & 0.003 & 0.258 & 0.02  & 90 \\
    3     & 2     & 3     & 0.004 & 0.037 & 0.02  & 150 \\
    4     & 2     & 4     & 0.007 & 0.197 & 0     & 50 \\
    5     & 3     & 6     & 0.004 & 0.018 & 0     & 50 \\
    6     & 4     & 5     & 0.004 & 0.037 & 0     & 130 \\
    7     & 5     & 6     & 0.002 & 0.140  & 0     & 50 \\
    \hline
    \end{tabular}
  \label{tab_line}
\end{table}

\begin{table}[htbp]\footnotesize
  \centering
  \caption{UPFC data.}
    \begin{tabular}{cccccc}
    \hline
    $T_{ij}^{sh,\max}$ & $T_{ij}^{se,\max}$ & $P_{ij}^{dc,\max}$ & $\Delta_{ij}^{P,U}$ & $\Delta_{ij}^{Q,se}$ & $\Delta_{ij}^{Q,sh}$ \\
    (MVA) & (MVA) & (MW)  & (MW)  & (MW)  & (MW) \\
    \hline
    100   & 100   & 100   & 200   & 200   & 200 \\
    \hline
    \end{tabular}
  \label{tab_upfc}
\end{table}%

Without loss of generality, the hourly wind power generation forecast error is assumed to follow a normal distribution $N(0,\sigma)$, and the standard deviation $\sigma$ is set as 20MW. This normal distribution assumption is an approximate and widely used one, which has been adopted in \cite{Wang2008,Wang2012,Pappala2009,Doherty2005}. Since the proposed UC model is independent of the distribution of wind power generation forecast error, other distributions may also be applied. 1000 scenarios of wind power generation are generated using LHS technique, each of which is assigned a probability that is one divided by the number of total generated scenarios, i.e., 0.001. Then these scenarios are reduced to 10 scenarios using the scenario reduction technique \cite{Morales2009}. It should be noted that the proposed UC models are independent of the scenario generation and reduction technique. As a result, other scenario generation and reduction methods can also be adopted. The reduced 10 scenarios are shown in Table \ref{tab_scen} in Appendix. More detailed scenario data can be found at \lj{\cite{Li2017}}.

In the following sections, \lj{the impacts of UPFC on wind power integration are analyzed from different aspects}, including wind power curtailment, load shedding, power flow, operating costs, unit status and voltage profile. \lj{Additionally, the approximate models are investigated}.

\subsection{Wind Power Curtailment and Load Shedding}
\begin{table}[htbp]\footnotesize
 \renewcommand{\arraystretch}{1.2}
  \centering
  \caption{Evaluation results of AC models.}
    \begin{tabular}{ccccc}
    \hline
    Model & EFC (\$)  & EWC (\$) & ELC (\$)  & ETC (\$)  \\
	\hline
    DM    & 106162.16  & 1469.22 & 1840.72  & 110219.76 \\
    NM    & 106084.34  & 1452.97 & 1646.59  & 109931.56 \\
    FSM   & 104461.65  & 1457.72 & 1801.01  & 108468.04  \\
    SSM   & 104309.68  & 383.77  & 137.23   & 105578.34  \\
    FSSM  & 102984.50  & 617.74  & 286.32   & 104636.22 \\
    \hline
    \end{tabular} \label{tab_ac_eva}
\end{table}

\begin{table}[htbp]\footnotesize
 \renewcommand{\arraystretch}{1.2}
  \centering
  \caption{Change rate of expected costs in AC models.}
    \begin{tabular}{ccccc}
    \hline
    Model & $CR^{EFC}$  & $CR^{EWC}$ & $CR^{ELC}$  & $CR^{ETC}$  \\
	\hline
    FSM   & -1.5\% & 0.3\%   & 9.4\%  & -1.3\% \\
    SSM   & -1.7\% & -73.6\% & -91.7\% & -4.0\% \\
    FSSM  & -2.9\% & -57.5\% & -82.6\% & -4.8\% \\
    \hline
    \end{tabular} \label{tab_ac_cr}
\end{table}

The expected costs are shown in Table \ref{tab_ac_eva}, while the change rates are listed in Table \ref{tab_ac_cr}. If UPFC is not allowed to be re-dispatched in the second stage (FSM), the impact on the expected wind power curtailment cost is negligible. Additionally, the expected load shedding cost increases compared with NM. Otherwise, dispatching UPFC in the second stage (SSM, FSSM) yields considerable reductions in the expected wind power curtailment and expected load shedding cost. It is especially evident in SSM, where the expected wind power curtailment cost and expected load shedding cost dramatically decrease by 73.5\% and 91.7\%, respectively. 

For further analysis, the wind power curtailment probability and loss of load probability are provided in Table \ref{tab_ac_wc_ls}. When observing FSM, one interesting finding is that employing UPFC may increase the probability of wind power curtailment and load shedding when the re-dispatch of UPFC is not allowed in the second stage. It implies that this traditional dispatch strategy, though reduces the total expected cost, may bring adverse impacts on wind power integration. In contrast, when dispatched in the second stage, UPFC significantly reduces the probability of wind power curtailment and load shedding by accommodating uncertain wind power generation. Another interesting observation is that SSM has lower probabilities of wind power curtailment and load shedding than FSSM, which coincides with the results of expected costs in Table \ref{tab_ac_eva}. The main reason behind this is that FSSM has a different UC schedule from SSM, yielding less expected total cost. In other words, FSSM achieves lower expected total cost at the expense of higher expected wind power curtailment cost and expected load shedding cost comparing with SSM.

\begin{table}[htbp]\footnotesize
 \renewcommand{\arraystretch}{1.2}
  \centering
  \caption{WPCP and LOLP of AC models.}
    \begin{tabular}{ccc}
    \hline
    Model &  WPCP & LOLP \\
	\hline
    DM    & 6.99\% & 4.85\% \\
    NM    & 6.87\% & 4.26\% \\
    FSM   & 7.51\% & 5.35\% \\
    SSM   & 2.12\% & 0.004\% \\
    FSSM  & 3.22\% & 0.72\% \\
    \hline
    \end{tabular} \label{tab_ac_wc_ls}
\end{table}

In summary, UPFC facilitates wind power integration when it is dispatched in the second stage, while it may not be able to adequately address uncertainty if it is only dispatched in the first stage. When considering only the expected wind power curtailment cost and expected load shedding cost, SSM is the best dispatch strategy.

\subsection{Power Flow}
\begin{figure}[htb]
\centering
\begin{tikzpicture}
\begin{axis}[
	% footnotesize,
	ybar,
	bar width=6pt,
	enlargelimits=0.15,
	ylabel={Wind Power Curtailment (MW)},
	symbolic x coords={NM, FSM, SSM, FSSM},
	xtick=data,
	nodes near coords,
	nodes near coords align={vertical},
]
	\addplot coordinates {(NM,21.6) (FSM,18.5) (SSM,0) (FSSM,0)};
	\footnotesize
\end{axis}
\end{tikzpicture}
\caption{Wind power curtailment in AC models.}
\label{fig_ac_wc}
\end{figure}
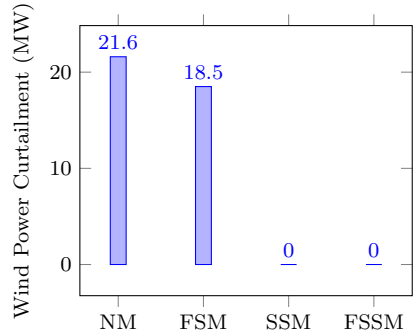

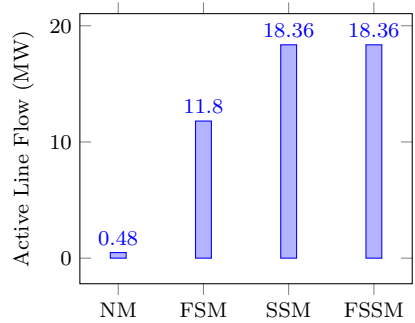
\begin{figure}[htb]
\centering
\begin{tikzpicture}
\begin{axis}[
	% footnotesize,
	ybar,
	bar width=6pt,
	enlargelimits=0.15,
	% legend style={at={(0.5,-0.15)},
	% anchor=north,legend columns=-1},
	ylabel={Active Line Flow (MW)},
	symbolic x coords={NM, FSM, SSM, FSSM},
	xtick=data,
	nodes near coords,
	nodes near coords align={vertical},
]
	\addplot coordinates {(NM,0.48) (FSM,11.80) (SSM,18.36) (FSSM,18.36)};
	% \addplot coordinates {(NM,130) (FSM,130) (SSM,130) (FSSM,130)};
	% \addplot coordinates {(NM,64) (FSM,60.47) (SSM,60.30) (FSSM,60.30)};
	% \addplot coordinates {(NM,0.48) (FSM,11.80) (SSM,18.36) (FSSM,18.36)};
	% \legend{Branch 4-5, Branch 1-4, Branch 4-2}
	\footnotesize
\end{axis}
\end{tikzpicture}
\caption{Active line flow through line 4-2.}
\label{fig_ac_pf}
\end{figure}

\lj{In order to further investigate the impacts of UPFC on power flow, hour 12 in scenario 78 is selected for analysis. At this hour, the active power flow through line 4-5 reaches the transmission capacity, which is 130MW, giving rise to a congestion, which causes 21.6 MW wind power curtailment in NM, as shown in Fig. \ref{fig_ac_wc}. However, less wind power generation is curtailed when UPFC is installed, since more active power is transferred through line 4-2 to bus 2, as depicted in Fig. \ref{fig_ac_pf}.} Particularly in SSM and FSSM, there is no wind power curtailment. This example shows that UPFC has the capability of regulating active power flow to reduce wind power curtailment.

\subsection{Operating Costs}

\begin{table*}[htbp]\footnotesize
  \centering
  \begin{threeparttable}
  \caption{Optimization results of AC models.}
    \label{tab_ac_opt}
    \begin{tabular}{cccccc}
    \hline
    \multirow{2}{*}{Model} & UC Cost & Expected  & Expected Wind Power & Expected Load  & Objective  \\
              &  (\$)     & Fuel Cost (\$)     & Curtailment Cost (\$) & Shedding Cost (\$)     & Cost\tnote{a} (\$) \\
    \hline
    DM    & 747.66 & 104649.90\tnote{b} & --    & --    & 105397.56 \\
    \lj{NOM}   & 1121.49 & 100072.66\tnote{b} & --    & --    & 101194.15 \\
    NM    & 747.66 & 106070.17 & 1345.07 & 1593.06 & 109755.96 \\
    FSM   & 747.66 & 104097.78 & 451.94 & 1863.59 & 107160.97 \\
    SSM   & 747.66 & 104384.35 & 123.19 & 33.68 & 105288.88 \\
    FSSM  & 747.66 & 103082.53 & 327.27 & 33.68 & 104191.14 \\
    \hline
    \end{tabular} 
    \begin{tablenotes}
        % \item [a] UC cost includes startup and shutdown cost of thermal unit.
        \item [a] Objective value of the optimization problem.
        % \item [b] Deterministic model (DM) without UPFC, based on wind power generation forecast.
        \item [b] The expected fuel cost of DM or NOM is based on wind power generation forecast.
    \end{tablenotes}
   \end{threeparttable}
\end{table*}

% \begin{table*}[htbp]\footnotesize
%   \centering
%   \begin{threeparttable}
%   \caption{Optimization results of AC models.}
%     \begin{tabular}{cccccc}
%     \hline
%     \multirow{2}{*}{Model} & UC Cost & Expected  & Expected Wind Power & Expected Load  & Objective  \\
%               &  (\$)     & Fuel Cost (\$)     & Curtailment Cost (\$) & Shedding Cost (\$)     & Cost\tnote{a} (\$) \\
%     \hline
%     DM    & 747.66 & 104649.90\tnote{b} & --    & --    & 105397.56 \\
%     NM    & 747.66 & 106070.17 & 1345.07 & 1593.06 & 109755.96 \\
%     FSM   & 747.66 & 104097.78 & 451.94 & 1863.59 & 107160.97 \\
%     SSM   & 747.66 & 104384.35 & 123.19 & 33.68 & 105288.88 \\
%     FSSM  & 747.66 & 103082.53 & 327.27 & 33.68 & 104191.14 \\
%     \hline
%     \end{tabular} \label{tab_ac_opt}
%     \begin{tablenotes}
%         % \item [a] UC cost includes startup and shutdown cost of thermal unit.
%         \item [a] Objective value of the optimization problem.
%         % \item [b] Deterministic model (DM) without UPFC, based on wind power generation forecast.
%         \item [b] The expected fuel cost of DM is based on wind power generation forecast.
%     \end{tablenotes}
%    \end{threeparttable}
% \end{table*}
The optimization results are provided in Table \ref{tab_ac_opt}. 
% The results indicate that UPFC has a negligible impact on UC cost, since the UC cost of most models are the same. However, UPFC does impact unit status, which will be discussed in the following section, where the UC schedule results are presented. 
\lj{Compared with DM, the utilization of UPFC in NOM helps to reduce fuel cost by changing the UC schedule.}
Comparing NM with DM, which is a deterministic model without UPFC, it clearly shows that the consideration of uncertainty increases the expected fuel cost if UPFC is not employed. On the contrary, with the help of UPFC (FSM, SSM, FSSM), the expected fuel cost falls to a lower level than DM. Meanwhile, the objective cost follows a similar trend. In particular, FSSM achieves the least expected fuel cost as well as objective cost since it is the most flexible model.

Different from the optimization results, in the evaluation results (see Table \ref{tab_ac_eva}), the expected fuel cost of DM is the highest due to its least robust UC schedule made with only wind power generation forecast. An interesting observation is no matter employed in which stage, UPFC has a beneficial effect on the expected fuel cost. It is particularly evident when UPFC is applied in both stages (FSSM), which leads to a 2.9\% reduction in the expected fuel cost. Therefore, it can be stated that when UPFC is applied, though UC cost remains unchanged, the expected fuel cost declines with wind power integration. Moreover, the best overall performance is achieved when UPFC is controllable in both stages (FSSM), resulting in the lowest expected total cost. \lj{One interesting finding is that when the UC schedule and UPFC setpoints are fixed as the solution to NOM, the economic dispatch problems may be infeasible in many scenarios. It indicates that the non-optimal method may result in high risk in operation.}

\subsection{Unit Status}
Though the UC costs in different models (DM, NM, FSM, SSM, FSSM) are the same, the unit status are different, as shown in Table \ref{tab_ac_uc}. The cheapest unit G1 is always committed regardless of model. Different from NM, G2 is committed at hour 8 and 9 when UPFC is employed. Additionally, when UPFC is used in the first stage (FSM, FSSM), G3 is committed for less hours. Especially in FSSM, G3 is only committed between hour 10 and 19. Consequently, the expected fuel cost of FSSM is reduced to achieve the smallest amount. It is understood that deploying UPFC leads to a more economic UC schedule.

% \begin{sidewaystable}\footnotesize
\begin{table*}[htbp]
  \footnotesize
  \centering
  \tabcolsep2pt
  \caption{\lj{UC results of AC models.}}
    \begin{tabular}{cccccccccccccccccccccccccc}
    \hline
          & Hour  & 1     & 2     & 3     & 4     & 5     & 6     & 7     & 8     & 9     & 10    & 11    & 12    & 13    & 14    & 15    & 16    & 17    & 18    & 19    & 20    & 21    & 22    & 23    & 24 \\
    \hline
    \multirow{5}{*}{G1} & DM  & 1     & 1     & 1     & 1     & 1     & 1     & 1     & 1     & 1     & 1     & 1     & 1     & 1     & 1     & 1     & 1     & 1     & 1     & 1     & 1     & 1     & 1     & 1     & 1 \\
    & NM    & 1     & 1     & 1     & 1     & 1     & 1     & 1     & 1     & 1     & 1     & 1     & 1     & 1     & 1     & 1     & 1     & 1     & 1     & 1     & 1     & 1     & 1     & 1     & 1 \\
    & FSM   & 1     & 1     & 1     & 1     & 1     & 1     & 1     & 1     & 1     & 1     & 1     & 1     & 1     & 1     & 1     & 1     & 1     & 1     & 1     & 1     & 1     & 1     & 1     & 1 \\
   	& SSM   & 1     & 1     & 1     & 1     & 1     & 1     & 1     & 1     & 1     & 1     & 1     & 1     & 1     & 1     & 1     & 1     & 1     & 1     & 1     & 1     & 1     & 1     & 1     & 1 \\
    & FSSM  & 1     & 1     & 1     & 1     & 1     & 1     & 1     & 1     & 1     & 1     & 1     & 1     & 1     & 1     & 1     & 1     & 1     & 1     & 1     & 1     & 1     & 1     & 1     & 1 \\
    \hline
    \multirow{5}[2]{*}{G2} & DM    & 1     & 0     & 0     & 0     & 0     & 0     & 0     & 0     & 0     & 1     & 1     & 1     & 1     & 1     & 1     & 1     & 1     & 1     & 1     & 1     & 1     & 1     & 1     & 1 \\
    & NM    & 1     & 0     & 0     & 0     & 0     & 0     & 0     & 0     & 0     & 1     & 1     & 1     & 1     & 1     & 1     & 1     & 1     & 1     & 1     & 1     & 1     & 1     & 1     & 1 \\
    & FSM   & 1     & 0     & 0     & 0     & 0     & 0     & 0     & 1     & 1     & 1     & 1     & 1     & 1     & 1     & 1     & 1     & 1     & 1     & 1     & 1     & 1     & 1     & 1     & 1 \\
    & SSM   & 1     & 0     & 0     & 0     & 0     & 0     & 0     & 1     & 1     & 1     & 1     & 1     & 1     & 1     & 1     & 1     & 1     & 1     & 1     & 1     & 1     & 1     & 1     & 1 \\
    & FSSM  & 1     & 0     & 0     & 0     & 0     & 0     & 0     & 1     & 1     & 1     & 1     & 1     & 1     & 1     & 1     & 1     & 1     & 1     & 1     & 1     & 1     & 1     & 1     & 1 \\
    \hline
    \multirow{5}[2]{*}{G3} & DM    & 0     & 1     & 1     & 1     & 1     & 1     & 1     & 1     & 1     & 1     & 1     & 1     & 1     & 1     & 1     & 1     & 1     & 1     & 1     & 1     & 1     & 0     & 0     & 0 \\
    & NM    & 0     & 1     & 1     & 1     & 1     & 1     & 1     & 1     & 1     & 1     & 1     & 1     & 1     & 1     & 1     & 1     & 1     & 1     & 1     & 1     & 1     & 1     & 0     & 0 \\
    & FSM   & 0     & 0     & 1     & 0     & 0     & 0     & 0     & 0     & 0     & 1     & 1     & 1     & 1     & 1     & 1     & 1     & 1     & 1     & 1     & 1     & 1     & 1     & 0     & 0 \\
    & SSM   & 0     & 1     & 1     & 1     & 1     & 1     & 1     & 0     & 0     & 1     & 1     & 1     & 1     & 1     & 1     & 1     & 1     & 1     & 1     & 0     & 0     & 1     & 0     & 0 \\
    & FSSM  & 0     & 0     & 0     & 0     & 0     & 0     & 0     & 0     & 0     & 1     & 1     & 1     & 1     & 1     & 1     & 1     & 1     & 1     & 1     & 0     & 0     & 0     & 0     & 0 \\
    \hline
    \end{tabular}%
  \label{tab_ac_uc}%
\end{table*}%
% \end{sidewaystable}

\subsection{Voltage}
\lj{Table. \ref{tab_ac_v} lists the evaluation results of the expected voltage magnitude at bus 4, which is connected with a wind farm. The results indicate that the average voltage magnitude is higher when UPFC is installed. Additionally, less fluctuation is experienced in SSM and FSSM than in FSM, since UPFC can be re-dispatched in the second stage according to different wind power generation scenarios.}

\begin{table}[htbp]\footnotesize
 \renewcommand{\arraystretch}{1.2}
  \centering
  \caption{Average voltage magnitude and variance of voltage magnitude.}
    \begin{tabular}{ccc}
    \hline
    Model &  Average Voltage Magnitude (p.u.) & Variance of Voltage Magnitude (p.u.) \\
	\hline
    NM    & 0.983 & 0.0002 \\
    FSM   & 0.993 & 0.0010 \\
    SSM   & 1.006 & 0.0005 \\
    FSSM  & 1.004 & 0.0004 \\
    \hline
    \end{tabular} \label{tab_ac_v}
\end{table}

\lj{In order to further investigate the impact of UPFC on voltage magnitude, the lower limit of voltage magnitude is increased from 0.95 p.u. to 0.98 p.u.}, while the upper limit remains at 1.05 p.u.. The optimization results are listed in Table \ref{tab_ac_opt_v}. The more rigid voltage limits lead to an increase in the objective cost (see Table \ref{tab_ac_opt}), while the impacts of UPFC are similar. Compared with NM, the expected wind power curtailment cost and load shedding cost are both reduced when UPFC is controllable in the second stage (SSM, FSSM). However, the expected load shedding cost increases in FSM, mainly because the power injections of UPFC are fixed in the second stage. Due to the reductions in the expected fuel cost and expected wind power curtailment cost, FSM still yields a lower objective cost than NM. One interesting observation is that FSSM outperforms SSM with respect to the expected wind power curtailment cost and expected load shedding cost, implying the impact of UC schedule becomes more significant.

\begin{table*}[htbp]\footnotesize
  \centering
  \caption{Optimization results of AC models with a more rigid voltage limit.}
    \begin{tabular}{cccccc}
    \hline
    \multirow{2}{*}{Model} & UC Cost & Expected  & Expected Wind Power & Expected Load  & Objective  \\
              &  (\$)     & Fuel Cost (\$)     & Curtailment Cost (\$) & Shedding Cost (\$)     & Cost\tnote{a} (\$) \\
    \hline
    NM    & 747.66 & 107229.50 & 1168.78 & 1343.58 & 110489.52 \\
    FSM   & 747.66 & 104142.27 & 455.41 & 1913.27 & 107258.61 \\
    SSM   & 1121.49 & 103696.46 & 321.29 & 476.45 & 105615.69 \\
    FSSM  & 747.66 & 103353.56 & 170.87 & 32.75 & 104304.84 \\
    \hline
    \end{tabular}
  \label{tab_ac_opt_v}
\end{table*}

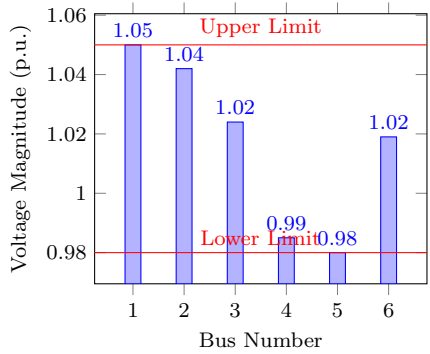
\begin{figure}[htb]
\centering
\begin{tikzpicture}
\begin{axis}[
	% footnotesize,
	ybar,
	bar width=6pt,
	enlargelimits=0.15,
	ylabel={Voltage Magnitude (p.u.)},
	% symbolic x coords={1, 2, 3, 4, 5, 6},
	xtick=data,
	xlabel={Bus Number},
	nodes near coords,
	nodes near coords align={vertical},
]
	\addplot coordinates {(1,1.05) (2,1.042) (3,1.024) (4,0.985) (5,0.980) (6,1.019)};
	\addplot[red,sharp plot,update limits=false] coordinates {(0,1.05) (7,1.05)}
	node[above] at (axis cs:3.5,1.05) {Upper Limit};
	\addplot[red,sharp plot,update limits=false] coordinates {(0,0.98) (7,0.98)}
	node[above] at (axis cs:3.5,0.98) {Lower Limit};
	\footnotesize
\end{axis}
\end{tikzpicture}
\caption{Voltage magnitude of all buses.}
\label{fig_case6_v}
\end{figure}

In some cases, the reactive power demand of the wind farm has to be curtailed due to the voltage magnitude limit. Consequently, wind power generation is also curtailed due to the constant power factor control strategy. Fig. \ref{fig_case6_v} illustrates the voltage magnitude of all the buses in such a situation, which occurs in NM in scenario 41 at hour 24. At this hour, G3 is off, whereas G1 and G2 are on. The voltage magnitude at bus 1 reaches its upper limit, while the reactive power output of G2 reaches its upper limit. Meanwhile, the voltage magnitude at bus 5 drops to its lower limit. In this situation, neither G1 nor G2 is able to provide sufficient reactive power to meet the demand of the wind farm. As a result, 1.8 MVar reactive power demand of the wind farm has to be curtailed, causing 6.2 MW curtailment of wind power generation.
% \begin{figure}[htb]
% \centering
% \includegraphics[width=2.5in]{case6v.pdf}
% \caption{Voltage magnitude of all the buses.}
% \label{fig_case6_v}
% \end{figure}

The results in this section demonstrate the capability of UPFC to improve voltage profile, and the impacts of UPFC when a more rigid voltage magnitude limit is imposed.

\subsection{DC and Mixed Models}
The optimization results of the DC models are listed in Table \ref{tab_dc_opt}. Similar to the AC models, UPFC can help to make a more economic UC schedule when used in the first stage, and reduces the expected wind power curtailment cost and expected load shedding cost when it is flexible in the second stage. The results of SSM and FSSM are the same, even though the setpoints of UPFC are different in the first stage. This is reasonable as there is no re-dispatch cost of UPFC. It also indicates that the flexibility of UPFC may not be fully exploited when only the DC constraints are considered.
\begin{table*}[htbp]\footnotesize
  \centering
  \caption{Optimization results of DC models.}
    \begin{tabular}{cccccc}
    \hline
    \multirow{2}{*}{Model} & UC Cost & Expected  & Expected Wind Power & Expected Load  & Objective  \\
              &  (\$)     & Fuel Cost (\$)     & Curtailment Cost (\$) & Shedding Cost (\$)     & Cost\tnote{a} (\$) \\
    \hline
    NM    & 747.66 & 106541.85 & 1540.03 & 603.05 & 109432.59 \\
    FSM   & 747.66 & 103059.14 & 409.32 & 2245.49 & 106461.61 \\
    SSM   & 747.66 & 102966.94 & 283.05 & 23.78 & 104021.43 \\
    FSSM  & 747.66 & 102966.94 & 283.05 & 23.78 & 104021.43 \\
    \hline
    \end{tabular}
  \label{tab_dc_opt}
\end{table*}

% \begin{table}[htbp]\footnotesize
%   \centering
%   \caption{Evaluation results of DC models.}
%     \begin{tabular}{ccccc}
%     \hline
%     Model & EFC (\$)  & EWC (\$) & ELC (\$)  & ETC (\$)  \\
% 	\hline
%     NM    & 106741.45 & 2178.61 & 807.54 & 110475.26 \\
%     FSM   & 103262.56 & 1260.46 & 2077.13 & 107347.81 \\
%     SSM   & 102918.32 & 897.73 & 232.36 & 104796.07 \\
%     FSSM  & 102918.32 & 897.73 & 232.36 & 104796.07 \\
%     \hline
%     \end{tabular}
%   \label{tab_dc_eva}
% \end{table}
The evaluation results of the DC models are shown in Table \ref{tab_dc_cr} and Table \ref{tab_dc_wc_ls}, which coincide with the optimization results. One interesting observation is that the expected load shedding cost in FSM apparently increases compared with that in NM, implying that the DC model shows the trend of change but may not be accurate.
\begin{table}[H]\footnotesize
 \renewcommand{\arraystretch}{1.2}
  \centering
  \caption{Evaluation results of DC models.}
    \begin{tabular}{ccccc}
    \hline
    Model & $CR^{EFC}$  & $CR^{EWC}$ & $CR^{ELC}$  & $CR^{ETC}$  \\
	\hline
    FSM   & -3.3\% & -42.1\% & 157.2\% & -2.8\% \\
    SSM   & -3.6\% & -58.8\% & -71.2\% & -5.1\% \\
    FSSM  & -3.6\% & -58.8\% & -71.2\% & -5.1\% \\
    \hline
    \end{tabular}
  \label{tab_dc_cr}
\end{table}

\begin{table}[H]\footnotesize
 \renewcommand{\arraystretch}{1.2}
  \centering
  \caption{WPCP and LOLP of DC models.}
    \begin{tabular}{ccc}
    \hline
    Model &  WPCP & LOLP \\
	\hline
    NM    & 9.60\% & 2.13\% \\
    FSM   & 6.34\% & 6.01\% \\
    SSM   & 4.30\% & 0.54\% \\
    FSSM  & 4.30\% & 0.54\% \\
    \hline
    \end{tabular} \label{tab_dc_wc_ls}
\end{table}

Table \ref{tab_mix_opt} and \ref{tab_mix_eva} present the optimization and evaluation results of the mixed models, respectively. The optimization problems of NM and SSM are infeasible, but they becomes feasible when the range of voltage magnitude limits is relaxed as 0.9 to 1.1, implying that in these cases the constraints associated with voltage cannot be ignored. When UPFC is employed in the first stage, the mixed models may be an approximate method due to the similar results to those of the AC models (see Table \ref{tab_ac_opt} and \ref{tab_ac_eva}). The WPCP and LOLP of the mixed models are shown in Table \ref{tab_mix_wc_ls}, which are also similar to the results of the AC models (see Table \ref{tab_ac_wc_ls}).
\begin{table*}[htbp]\footnotesize
  \centering
  \caption{Optimization results of mixed models.}
    \begin{tabular}{cccccc}
    \hline
    \multirow{2}{*}{Model} & UC Cost & Expected  & Expected Wind Power & Expected Load  & Objective  \\
              &  (\$)     & Fuel Cost (\$)     & Curtailment Cost (\$) & Shedding Cost (\$)     & Cost\tnote{a} (\$) \\
    \hline
    NM    & --    & --    & --    & --    & -- \\
    FSM   & 747.66 & 104097.78 & 451.94 & 1863.59 & 107160.97 \\
    SSM   & --    & --    & --    & --    & -- \\
    FSSM  & 747.66 & 103282.54 & 152.39 & 33.68 & 104216.27 \\
    \hline
    \end{tabular}%
  \label{tab_mix_opt}%
\end{table*}%

\begin{table}[H]\footnotesize
  \centering
  \caption{Evaluation results of mixed models.}
    \begin{tabular}{ccccc}
    \hline
    Model & EFC (\$)  & EWC (\$) & ELC (\$)  & ETC (\$)  \\
	\hline
    NM    & --    & --    & --    & -- \\
    FSM   & 104461.65 & 1457.72 & 1801.01 & 108468.04 \\
    SSM   & --    & --    & --    & -- \\
    FSSM  & 103314.1 & 652.45 & 281.41 & 104995.62 \\
    \hline
    \end{tabular}
  \label{tab_mix_eva}
\end{table}

\begin{table}[H]\footnotesize
 \renewcommand{\arraystretch}{1.2}
  \centering
  \caption{WPCP and LOLP of mixed models.}
    \begin{tabular}{ccc}
    \hline
    Model &  WPCP & LOLP \\
	\hline
    NM    & -- & -- \\
    FSM   & 7.51\% & 5.35\% \\
    SSM   & -- & -- \\
    FSSM  & 3.28\% & 0.63\% \\
    \hline
    \end{tabular} \label{tab_mix_wc_ls}
\end{table}

The computational time of different models is listed in Table \ref{tab_cal_time}. The computational time of DC models declines dramatically due to much less complexity. Meanwhile, the mixed models have longer computational time than the DC models, but still much shorter than the AC models.
% \begin{table}[H]\footnotesize
%   \centering
%   \caption{Computational time of different models.}
%     \begin{tabular}{cccc}
%     \hline
%          & AC (s) & DC (s) & Mixed (s) \\
%     \hline
%     NM   & 86.23  & 2.04  & -- \\
%     FSM  & 123.40 & 2.29  & 12.34 \\
%     SSM  & 190.32 & 1.74  & -- \\
%     FSSM & 170.31 & 1.94  & 11.91 \\
%     \hline
%     \end{tabular}%
%   \label{tab_cal_time}%
% \end{table}%

\begin{table}[htbp]\footnotesize
  \centering
  \caption{Computational time of different models.}
    \begin{tabular}{cccc}
    \hline
         & AC (s) & DC (s) & Mixed (s) \\
    \hline
    NM   & 86.23\tnote{*}  & 2.04  & -- \\
    FSM  & 123.40 & 2.29  & 12.34 \\
    SSM  & 190.32 & 1.74  & -- \\
    FSSM & 170.31 & 1.94  & 11.91 \\
    \hline
    \end{tabular}%
  \label{tab_cal_time}%
      \begin{tablenotes}
 \item  *The DICOPT solver, which is based on outer approximation algorithm, is used to solve the original mixed-integer non-linear problems (AC). Then the computational times are compared with the approximate DC and mixed models (Mixed).  \end{tablenotes}
\end{table}%

\begin{figure}[H]
\centering
\begin{tikzpicture}
	\begin{axis}[
		xlabel=Capacity of UPFC,
		ylabel=Objective Cost (\$),
		% legend style={at={(0.5,-0.20)},
		% anchor=north,legend columns=-1},
	]
	\addplot coordinates {
	(50,105222.06) (100,104191.14) (150,104087.10) (200,104050.75)
	};
	\addplot coordinates {
	(50,105641.53) (100,104021.43) (150,103431.70) (200,103355.08)
	};
	\addplot coordinates {
	(50,105222.06) (100,104216.27) (150,104103.85) (200,104080.84)
	};
	\legend{AC, DC, Mixed}
	\footnotesize
	\end{axis}
\end{tikzpicture}
\caption{Optimization results with increasing UPFC capacity.}
\label{fig_compare_obj}
\end{figure}
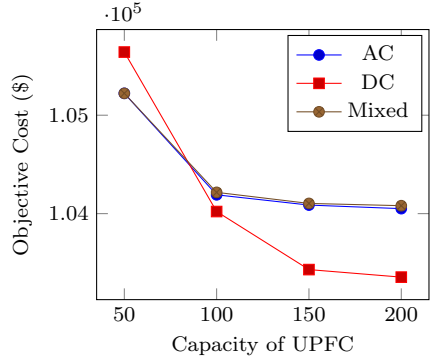

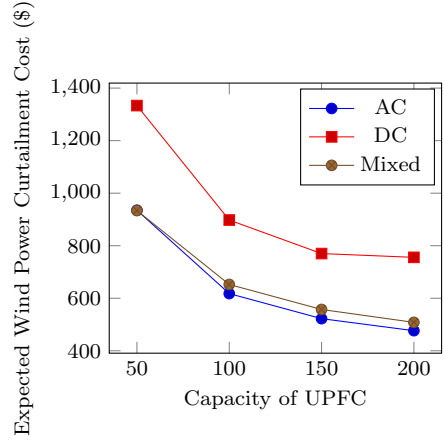
\begin{figure}[H]
\centering
\begin{tikzpicture}
	\begin{axis}[
		xlabel=Capacity of UPFC,
		ylabel=Expected Wind Power Curtailment Cost (\$),
		% legend style={at={(0.5,-0.20)},
		% anchor=north,legend columns=-1},
	]
	\addplot coordinates {
	(50,934.94) (100,617.74) (150,522.31) (200,476.56)
	};
	\addplot coordinates {
	(50,1333.41) (100,897.73) (150,770.09) (200,755.72)
	};
	\addplot coordinates {
	(50,933.79) (100,652.45) (150,557.26) (200,508.19)
	};
	\legend{AC, DC, Mixed}
	\footnotesize
	\end{axis}
\end{tikzpicture}
\caption{Expected wind power curtailment cost with increasing UPFC capacity.}
\label{fig_compare_wc}
\end{figure}

\begin{figure}[H]
\centering
\begin{tikzpicture}
	\begin{axis}[
		xlabel=Capacity of UPFC,
		ylabel=Expected Load Shedding Cost (\$),
		% legend style={at={(0.5,-0.20)},
		% anchor=north,legend columns=-1},
	]
	\addplot coordinates {
	(50,502.31) (100,286.32) (150,206.59) (200,186.58)
	};
	\addplot coordinates {
	(50,384.90) (100,232.36) (150,187.32) (200,186.01)
	};
	\addplot coordinates {
	(50,499.90) (100,281.41) (150,201.93) (200,185.58)
	};
	\legend{AC, DC, Mixed}
	\footnotesize
	\end{axis}
\end{tikzpicture}
\caption{Expected load shedding cost with increasing UPFC capacity.}
\label{fig_compare_ls}
\end{figure}
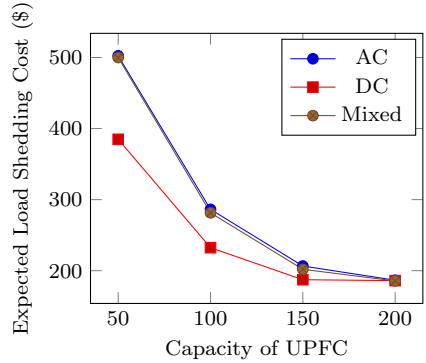

In order to demonstrate the trend of solutions with changing parameters, \lj{FSSM is investigated}, which is the best overall performing approach, in AC, DC, and mixed formulations with different capacities of UPFC. When the capacity of UPFC increases, the objective cost decreases as shown in Fig. \ref{fig_compare_obj}, and the deviation between the AC and DC model becomes larger, while the results of the mixed model remains close to that of the AC model. As depicted in Fig. \ref{fig_compare_wc} and \ref{fig_compare_ls}, the evaluation results show that increasing the capacity of UPFC leads to lower expected wind power curtailment cost and expected load shedding cost. The interesting point is that the impact on these costs declines as the capacity of UPFC increases, implying that the benefit of UPFC gets smaller when its capacity further increases. In addition, the expected wind power curtailment cost of the DC model is always higher than that of the AC model, indicating that using the DC model may overestimate the expectation of wind power curtailment.

Therefore, it may be implied that compared with the AC model, the DC model has a much lower computational burden, but the objective cost may be inaccurate, especially when the capacity of UPFC is large. Meanwhile, the mixed model obtains relatively accurate results with higher computational efficiency. Thus, the mixed model may be an alternative way to solve the problem.

\section{Concluding Remarks}\label{sec_conclusion}
\lj{UPFC has the underlying capability to facilitate wind power integration in UC, which has not been fully exploited. In this paper, a comprehensive evaluation process based on a two-stage stochastic UC model with UPFC and wind power integration is proposed. Different ways of utilizing UPFC are compared to analyze the impacts and benefits of UPFC, for bettering understanding the role of UPFC in accommodating uncertain wind power generation. Some conclusions are drawn as follows.}

\lj{(1) When set up only in the first stage (FSM), UPFC has negligible impacts on the UC cost, but reduces the fuel cost by affecting unit status, leading to a more economic UC schedule. However, it may has adverse impacts on wind power integration, causing higher probabilities of wind power curtailment and load shedding due to the inflexible control in the second stage.}

\lj{(2) When dispatched only in the second stage (SSM), UPFC yields a considerable decrement of expected wind power curtailment cost as well as expected load shedding cost, due to its capability of regulating power flow. In fact, if only the expected wind power curtailment cost and expected load shedding cost are considered, SSM performs better than the other models within relatively relaxed voltage limits. Moreover, voltage profile is improved and the expected total cost is reduced even with more rigid voltage magnitude limits.}

\lj{(3) Reaching the lowest expected total cost, FSSM achieves the best performance for facilitating wind power integration in UC. In FSSM, the flexibility of UPFC is fully exploited so that all types of costs are reduced. However, within relatively relaxed voltage limits, FSSM sacrifices the savings in the expected wind power curtailment cost and expected load shedding cost for less expected fuel cost. As a result, power system operators may choose to employ UPFC in different stages according to different purposes. Further, the formulation is easy to be extended to model other FACTS devices.}

\lj{(4) Compared with the AC model, the DC model is incapable of taking full advantage of UPFC. In contrast, the mixed model may obtain similar results to the AC model with much less computational complexity. However, the case studies also show that, the mixed model is not applicable to NM and SSM, since AC power flow constraints are not satisfied with the UC schedules derived from the DC models.}

\lj{(5) The benefits of UPFC become smaller when the capacity of UPFC  increases to a certain extent. That implies transmission congestion is not the bottleneck of integrating wind generation any longer, and other types of controllable resources other than transmission flexibility are required.}

Technically, the AC problems may be solved using interior point algorithm \cite{Zhang2001a,Zhang2001b}, intelligent algorithm \cite{Panuganti2013,Leung2011,Bhasaputra2002,Chung2001}, or decomposition technique \cite{Wang2008,Shaoyun1998}. While this work serves as a starting point to solve the problems for a real-world large-scale system, the future work will focus on developing more efficient and practical algorithms to solve the problems.

\lj{Some emerging techniques, such as energy storage, also provide flexibility to accommodates the uncertainty of wind power generation. Energy storage systems accommodate uncertainty by re-balancing power generation and demand. They can either store surplus wind generation or release electricity when wind generation is insufficient. Such flexibility, however, has to be transfered through transmission network if generation and demand are not in the same place. In this context, transmission flexibility is required to eliminate possible congestions. Though FACTS device itself does not generate active power, it is specifically used for alleviating congestions by controlling power flows. In this sense, FACTS devices cannot replace energy storage, but rather serve as a supplement for fully exploiting power system flexibility. The coordination of FACTS devices and energy storage is well worth investigating for future power systems.}

% \section*{Acknowledgments}
% This work was supported in part by the Foundation for Innovative Research Groups of the National Natural Science Foundation of China (51621065), in part by the China State Grid Corp Science and Technology Project (SGSXDKY-DWKJ2015-001), and in part by the State Key Development Program of Basic Research of China (2013CB228201).

\clearpage
\section*{Appendix}
\begin{table}[H]\footnotesize
  \centering
  \caption{Generator data.}
    \begin{tabular}{cccccccccc}
    \hline
    Unit  & Bus  & Pmax & Pmin & Qmax & Qmin & Ini.  & Min  & Min & Ramp \\
          & No.  & (MW) & (MW) & (MVar) & (MVar) & State & Off & On & (MW/h) \\
          &   &  &  & &  & (h) & (h) & (h) &  \\
    \hline
    G1    & 1     & 220   & 90    & 200   & -80   & 4     & -4    & 4     & 50 \\
    G2    & 2     & 100   & 10    & 70    & -40   & 2     & -3    & 2     & 40 \\
    G3    & 6     & 20    & 10    & 70    & -40   & -1    & -1    & 1     & 15 \\
    \hline
    \end{tabular}%
  \label{tab_gen_data}%
\end{table}%

\begin{table}[H]\footnotesize
  \centering
  \caption{Generator operating cost data.}
    \begin{tabular}{cccccc}
    \hline
    Unit & \multicolumn{3}{c}{Fuel Consumption Function} & Startup Fuel & Fuel Price \\
         & a (MBTU) & b (MBTU/MWh) & c (MBTU/$\text{MW}^2$h) & (MBTU) & (\$/MBTU) \\
    \hline
    G1    & 176.9 & 13.5  & 0.0004 & 100   & 1.2469 \\
    G2    & 129.9 & 32.6  & 0.001 & 300   & 1.2461 \\
    G3    & 137.4 & 17.6  & 0.005 & 0     & 1.2462 \\
    \hline
    \end{tabular}%
  \label{tab_operation_data}%
\end{table}%

\begin{table}[H]\footnotesize
  \tabcolsep2pt
  \centering
  \caption{Hourly load and wind power generation forecast.}
    \begin{tabular}{cccccccc}
    \hline
    Hour  & Pload (MW) & Qload (MW) & Pwind (MW) & Hour  & Pload (MW) & Qload (MW) & Pwind (MW) \\
    \hline
    1     & 219.19 & 50.4  & 44    & 13    & 326.18 & 69.6  & 84 \\
    2     & 235.35 & 47.4  & 70.2  & 14    & 323.6 & 70    & 80 \\
    3     & 234.67 & 45.6  & 76    & 15    & 326.86 & 71.6  & 78 \\
    4     & 236.73 & 44.5  & 82    & 16    & 287.79 & 73.5  & 32 \\
    5     & 239.06 & 44.6  & 84    & 17    & 260   & 73.6  & 4 \\
    6     & 244.48 & 46.1  & 84    & 18    & 246.74 & 70.9  & 8 \\
    7     & 273.39 & 49.9  & 100   & 19    & 255.97 & 70.7  & 10 \\
    8     & 290.4 & 51.1  & 100   & 20    & 237.35 & 68.2  & 5 \\
    9     & 283.56 & 53.7  & 78    & 21    & 243.31 & 68.2  & 6 \\
    10    & 281.2 & 59.5  & 64    & 22    & 283.14 & 66.9  & 56 \\
    11    & 328.61 & 65.7  & 100   & 23    & 283.05 & 56.3  & 82 \\
    12    & 328.1 & 67.9  & 92    & 24    & 248.75 & 56.2  & 52 \\
    \hline
    \end{tabular}%
  \label{tab_load}%
\end{table}%

\begin{table}[H]\footnotesize
  \centering
  \caption{Reduced wind power generation scenarios.}
    \begin{tabular}{ccccccccccc}
    \hline
    Hour & \multicolumn{10}{c}{Wind Power Generation Scenarios} \\
          & 1     & 2     & 3     & 4     & 5     & 6     & 7     & 8     & 9     & 10 \\
    \hline
    1     & 79.2  & 39.8  & 40.4  & 35.2  & 54.1  & 34.4  & 81.7  & 37.5  & 58.4  & 50.4 \\
    2     & 79.7  & 71.7  & 55.5  & 99.6  & 64.6  & 67    & 98.9  & 80    & 70.9  & 68.8 \\
    3     & 66    & 66    & 33.9  & 68.6  & 92.5  & 34.3  & 52.2  & 92.1  & 60.2  & 69.4 \\
    4     & 65.8  & 87.3  & 74.6  & 94    & 77.7  & 106.6 & 97.4  & 69.9  & 86.4  & 73.4 \\
    5     & 93    & 78.8  & 93.4  & 103.2 & 90.1  & 64.3  & 78.6  & 84.5  & 116.2 & 65.5 \\
    6     & 88.3  & 96.7  & 111   & 54    & 69.8  & 85.9  & 93.1  & 111.2 & 99.4  & 94.2 \\
    7     & 111.2 & 111.1 & 110.2 & 101.5 & 92.5  & 112.4 & 118.1 & 118.8 & 103.7 & 99 \\
    8     & 126.4 & 59.5  & 64.9  & 94.6  & 77.9  & 113.2 & 92.5  & 98.1  & 124.1 & 97.4 \\
    9     & 77.1  & 96.4  & 59.8  & 77    & 98.2  & 50.5  & 85.3  & 89.8  & 65.1  & 55.3 \\
    10    & 81    & 56.3  & 74.3  & 40.6  & 87.5  & 46.7  & 58    & 75.7  & 50.7  & 30.5 \\
    11    & 116.2 & 135.2 & 54.2  & 117.8 & 124.9 & 65.9  & 120.2 & 119.2 & 77.7  & 34.2 \\
    12    & 90.9  & 85.2  & 63.8  & 99.8  & 112.9 & 78.7  & 63.5  & 68.7  & 124   & 101.6 \\
    13    & 89.1  & 87.6  & 89.9  & 81.8  & 92.6  & 104.2 & 75.7  & 59.1  & 80.7  & 136.8 \\
    14    & 75.7  & 65    & 54.4  & 88.8  & 40.7  & 68.1  & 101.7 & 43.3  & 84.3  & 75.3 \\
    15    & 63.1  & 81.6  & 64.9  & 99.5  & 100.7 & 82.1  & 94.3  & 65.9  & 52.5  & 82.6 \\
    16    & 24.6  & 0     & 39.2  & 47.4  & 11.1  & 14.3  & 44    & 43.5  & 45.2  & 0 \\
    17    & 9.7   & 14.2  & 0     & 22.9  & 0     & 0     & 25.1  & 15    & 2     & 0 \\
    18    & 4.8   & 0     & 0     & 24.6  & 3     & 12.6  & 13.5  & 2.8   & 28    & 37.1 \\
    19    & 0     & 0     & 0     & 0     & 29.8  & 40    & 32.5  & 0     & 0     & 13.2 \\
    20    & 0     & 20.5  & 0     & 1.8   & 9.5   & 28.6  & 20.2  & 0     & 0     & 8.9 \\
    21    & 0     & 0     & 11.7  & 70.6  & 0     & 29.4  & 0     & 15.8  & 2.9   & 11.1 \\
    22    & 76.1  & 19.3  & 69.2  & 76.4  & 41.5  & 71.5  & 22.3  & 34.3  & 65.4  & 22.3 \\
    23    & 78.3  & 92.6  & 84    & 111.7 & 112.3 & 102.5 & 72.4  & 107.5 & 77.8  & 79.7 \\
    24    & 56.4  & 30    & 34.2  & 53.9  & 39.3  & 25.4  & 47.3  & 33.3  & 18.1  & 37.1 \\
    \hline
    \end{tabular}
  \label{tab_scen}
\end{table}

\clearpage
\section*{References}

% \bibliography{library}
\bibliography{../../../library}

\end{document}